\documentclass[12pt,leqno]{article}
\usepackage{graphicx}
\font\tendb=msbm10 at 12pt \font\sevendb=msbm10 at 9pt
\font\fivedb=msbm10 at 7pt
\newfam\dbfam
\textfont\dbfam=\tendb \scriptfont\dbfam=\sevendb
\scriptscriptfont\dbfam=\fivedb
\def\db{\fam\dbfam\tendb}
\font\eufm=eufm10\font\eufms=eufm10\font\eufmss=eufm10\newfam\eufam
\textfont\eufam=\eufm\scriptfont\eufam=\eufms\scriptscriptfont\eufam=\eufmss

\font\tendbb=msbm10 at 12pt \font\sevendbb=msbm7 at 9pt
\font\fivedbb=msbm5 at 6pt
\newfam\dbbfam
\textfont\dbbfam=\tendbb \scriptfont\dbbfam=\sevendbb
\scriptscriptfont\dbbfam=\fivedbb

 \def \Z {{\db Z}}
 \def \R {\hbox{\db R}}

 \def \S {S^{3}}

\def \fin {\hfill \framebox(7,7) }
 
\font\tenMmm=eusm10 at 12pt
\newfam\Mmmfam

\textfont\Mmmfam=\tenMmm

\def\illu #1 by #2 (#3){
  \vbox to #2{
    \hrule width #1 height 0pt depth 0pt
    \vfill
    \special{illustration #3} 
    }
  }

\begin{document}

\null \vspace{1cm}

\begin{center}
{\large {\bf Skein Algebras of the solid torus and  symmetric spatial graphs}}\\
 Nafaa Chbili\footnote{ Supported by a  fellowship from the Japan Society for the
Promotion of Science (JSPS), and by a Grant-in-Aid for JSPS
fellows 03020. The author would
like to thank the JSPS for its support.}\\
\begin{footnotesize}
Department of Mathematics\\
 Tokyo Institute of Technology\\
 Oh-okayama Meguro Tokyo 152-8551 JAPAN\\
 E-mail: chbili@math.titech.ac.jp\\
  \end{footnotesize}
\end{center}

\begin{abstract}
We use the topological invariant of spatial graphs introduced by S. Yamada to find
necessary conditions for a spatial graph to be periodic with  a prime period.
The proof of the main result  is based on computing the Yamada skein algebra of the solid
 torus
then proving that   this algebra injects into the Kauffman bracket skein
 algebra of the solid torus.
\end{abstract}

\begin{center}
{\sc I- Introduction}
\end{center}
Throughout this paper a  graph  is a finite one-dimensional
CW-complex. We assume that vertices of  our graphs  have valency
greater than or equal to 3. If $G$ is a graph, we denote by
${\mathcal V }(G)$ the set of vertices and by ${\mathcal E}(G)$
the set of edges of $G$. An embedding of such a graph into the
three-sphere is called a spatial graph. Let $p \geq 2$ be an
integer. A spatial graph $\tilde G$ is said to be $p$-periodic if
there exists an orientation preserving diffeomorphism $h$ of
$(\S,\tilde G)$ such that $h$ is of order $p$ and  the set of
fixed points of $h$ is a circle that does not intersects $\tilde
G$. If $\tilde G$ is a periodic spatial graph, we will denote the
quotient spatial graph by $\underline {\tilde G}$. \\
By the positive solution of the Smith conjecture \cite{BM}, the
action defined by $h$ on the three-sphere  is topologically
conjugate to an orthogonal action. In other words, if we identify
the three-sphere with $\R^3 \cup \infty$, then we may assume that
$h$ is a rotation by a $2\pi/p$ angle around the standard $z$-axis
in the Euclidean space $\R^3$. Hence, a $p$-periodic spatial graph
may be represented by a diagram in $\R^2$ which is invariant by a planar rotation. \\
Let $\tilde G$ be a $p$-periodic spatial graph. We define the
wrapping number of $\tilde G$ to be the minimum number of
intersection points of $\tilde G$ with  $P$, where $P$ runs
over all half-planes bounded by the axis of the rotation \cite{Ma}.\\
The periodicity of knots and links have been subject to extensive
literature. Both  the new and the classical invariants have been
used successfully   to study the periodicity of knots and links.
In particular, several criteria for periodicity of links have been
introduced in terms of the invariants of skein type. For instance,
using the Jones polynomial (or, equivalently, the Kauffman
bracket), Murasugi \cite{Mu} and  Traczyk \cite{Tr} introduced
necessary conditions for a link to be periodic with  a prime
period. These criteria have been extended to the skein polynomial
(HOMFLYPT) by Przytycki \cite {Pr1}.\\
In \cite {Ya1}, Yamada introduced a topological invariant of
 spatial graphs which takes values in
 ${\mathcal R}=\Z[A^{\pm 1}, d^{-1}]$, where
$d=-A^{2}-A^{-2}$. This invariant, denoted here by $Y$, can be
defined recursively on planar diagrams of spatial graphs. When
restricted to knots, this invariant is essentially  the Jones
polynomial  of some cable of the knot. The goal of this paper is
the use of the Yamada polynomial to find necessary conditions for
a spatial graph to be periodic. In the special cases corresponding
to periodic graphs with wrapping number one or two, Marui
\cite{Ma} introduced  necessary conditions, in terms of the Yamada
polynomial, for a spatial graph to be periodic with wrapping
number one or two. The main result in our paper is given by the
following theorem, where
the congruences hold in the ring $\Z[A^{\pm 1},d^{-1}]$.\\

{\bf Theorem 1.1} {\sl Let $p$  be a prime  and $\tilde G$ a
spatial graph. If $\tilde G$ is $p-$periodic, then}\\
\begin{tabular}{ll}
(1)& $Y(\tilde G)(A)\equiv (Y(\tilde {\underline G})(A))^p$
 modulo $p, d^{p}-d$.\\
(2)& $Y(\tilde G)(A)\equiv Y(\tilde G)(A^{-1})$ modulo $p,
A^{2p}-1$.
\end{tabular}
\\

 {\bf Application.}  Let $P$ be the Peterson graph and  $\tilde P$ be
the embedding of $P$ given by  Figure 1. Marui \cite{Ma}, proved
that this spatial graph cannot be 5-periodic with wrapping number
two.

\begin{center}
\includegraphics[width=4cm,height=5cm]{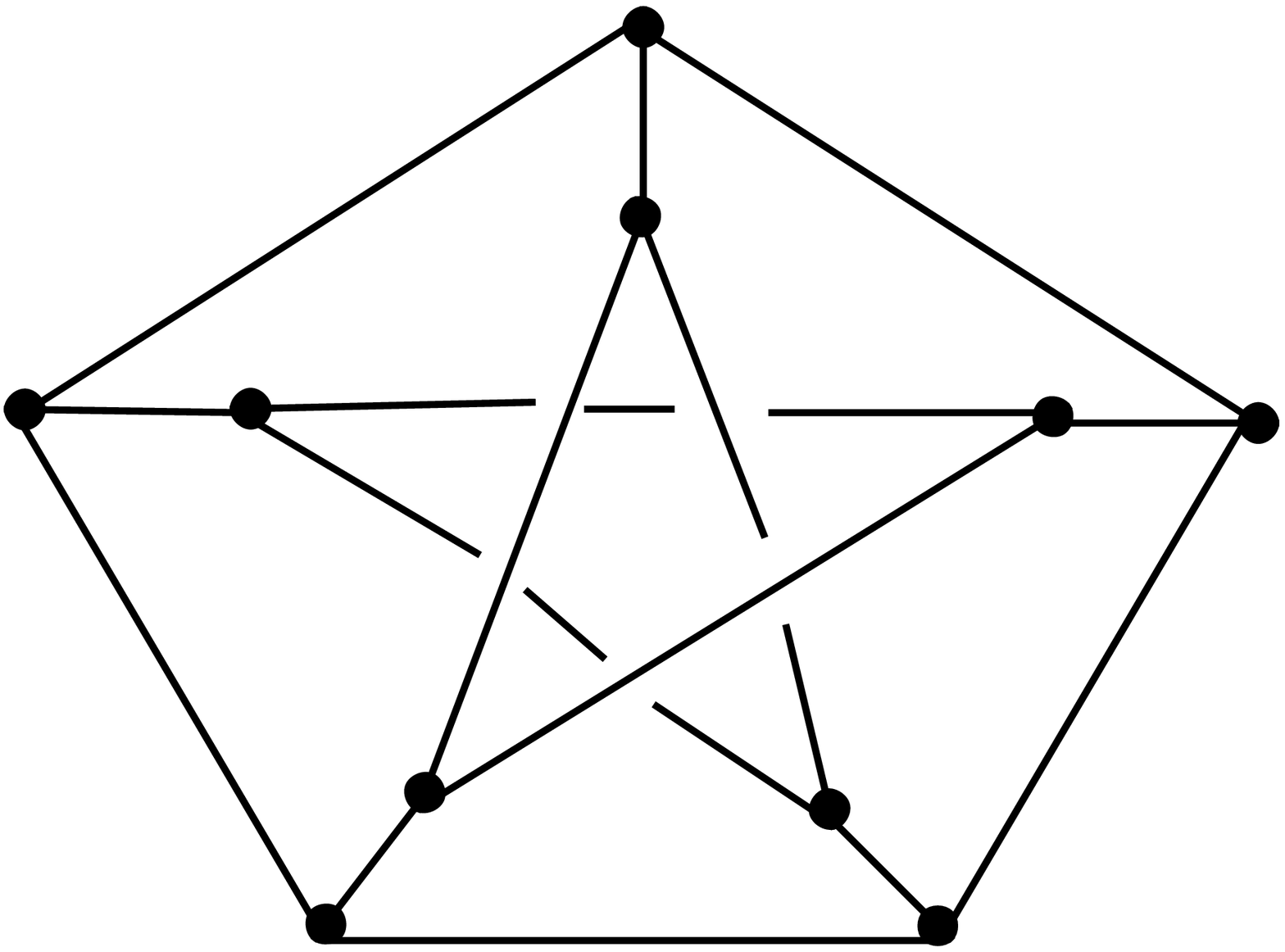}
\end{center}

\begin{center} {\sc Figure 1} \end{center}
If we apply Theorem 1.1, (2) we can prove that $\tilde P$ is not
5-periodic with any wrapping number. From the computation in \cite
{Ma}, we have the following congruence modulo 5, $A^{10}-1$:
 $$Y(\tilde P)(A) \equiv 2(4A^{6}+4+A^{4}
 +A^{8}).$$
Obviously this polynomial does not satisfy condition (2) of
Theorem 1.1. Hence,
the spatial graph $\tilde P$ is not 5-periodic.\\
It is worth mentioning here that the finite group $\Z/5\Z$ acts
freely  on the  abstract graph $P$.\\
{\bf Remark.}
 Let  $G$ be a graph and let
 $\tilde G$ be  a spatial embedding of $G$ which is $p-$periodic. Then, the rotation
  induces
actions of the finite group $\Z/p\Z$ on  ${\mathcal E}(G)$ and
${\mathcal V}(G)$. Obviously,  the action of $\Z/p\Z$ on
${\mathcal V}(G)$ is free. Moreover, vertices that belong to the
same orbit have the same valency.
 Hence, we can get some restriction on
the possible periods for a spatial graph $\tilde G$ from the
properties of the original abstract graph. In particular, from the
number and the valency of vertices. For instance, we can easily
see
 that for all integers $p$ an embedding of a $\theta $-curve
cannot be $p$-periodic. Also, an embedding of the complete graph
$K_n$ cannot be 2-periodic.
\begin{center} {\sc II- The Yamada Polynomial} \end{center}
 A ribbon graph is an oriented compact
surface with boundary that retracts by deformation on a graph. A
spatial ribbon graph is an embedding of a ribbon graph into the
three-sphere. It is well known that the study of spatial ribbon
graphs up to isotopy is equivalent to the study of planar graph
diagrams up to the extended Reidemeister moves \cite{Ya1}.\\
 In \cite{Ya1}, S. Yamada introduced an invariant $R$ of regular
isotopy of spatial graphs. This invariant takes its values in the
ring $\Z[A^{\pm 1}]$ and may be defined recursively on diagrams of
spatial graphs.  A similar invariant of trivalent graphs, with
good weight associated with the set of edges, was also introduced
by Yamada \cite{Ya2}. This invariant was extended by Yokota
\cite{Yo} using the linear skein theory introduced by Lickorish
\cite{Li}. In the present  paper, we find it more convenient to
slightly change the recursive formulas introduced by Yamada.
Namely, we define an invariant $Y$ of spatial graphs recursively
by the initialization $Y(\O)=1$  and  the four relations in Figure
2. Notice that the vertical dots in our figures  mean an arbitrary
number of edges. It is also worth mentioning that the following
identities hold for
diagrams which are identical except in a small disk where they look as indicated below.\\
\null

               \begin{picture}(0,0)

               \put(135,0){\line(-1,1){13}}
               \put(119,16){\line(-1,1){13}}
               \put(105,0){\line(1,1){30}}
               \put(160,10){$= A^{4}$$Y($}
               \put(80,10){$Y($}
               \put(150,10){$)$}
               \put(205,15){\oval(28,20)[r]}
               \put(235,15){\oval(28,20)[l]}
               \put(234,10){$)$$ +A^{-4}$$Y($}
               \put(300,25){\oval(30,20)[b]}
               \put(300,0){\oval(30,20)[t]}
               \put(324,10){$)$$ -d$ $Y($}
               \put(384,0){\line(-1,1){30}}
               \put(356,0){\line(1,1){30}}
               \put(387,10){$)$}

\end{picture}
\\

\begin{picture}(0,0)
                \put(120,10){\line(-1,1){13}}
                \put(110,12){.}
                \put(110,10){.}
                \put(110,8){.}

               \put(120,10){\line(-1,-1){13}}
               \put(120,10){\line(1,0){15}}
               \put(135,10){\line(1,1){13}}
               \put(135,10){\line(1,-1){13}}
               \put(145,12){.}
                \put(145,10){.}
                \put(145,8){.}
               \put(80,10){$Y($}
               \put(150,10){$)$}
               \put(160,10){$= Y($}
                \put(200,10){\line(-1,1){13}}
               \put(200,10){\line(-1,-1){13}}
               \put(190,12){.}
                \put(190,10){.}
                \put(190,8){.}
               \put(200,10){\line(1,1){13}}
               \put(200,10){\line(1,-1){13}}
               \put(210,12){.}
                \put(210,10){.}
                \put(210,8){.}
                \put(215,10){$)-d^{-1}Y($}
                \put(275,10){\line(-1,1){13}}
               \put(275,10){\line(-1,-1){13}}
               \put(265,12){.}
                \put(265,10){.}
                \put(265,8){.}
               \put(280,10){\line(1,1){13}}
               \put(280,10){\line(1,-1){13}}
               \put(290,12){.}
                \put(290,10){.}
                \put(290,8){.}
                \put(295,10){$)$}
\end{picture}
\\

\begin{picture}(0,0)
\put(120,10){\line(-1,1){13}}
                \put(110,12){.}
                \put(110,10){.}
                \put(110,8){.}
               \put(120,10){\line(-1,-1){13}}
               \put(130,10){\circle{20}}
                \put(80,10){$Y($}
               \put(150,10){$)$}
\put(160,10){$= (d-d^{-1})Y($} \put(245,10){\line(-1,1){13}}
                \put(235,12){.}
                \put(235,10){.}
                \put(235,8){.}
               \put(245,10){\line(-1,-1){13}}
               \put(250,10){$)$}
\end{picture}
\\

\begin{picture}(0,0)

\put(80,10){$Y( \;\;D \bigsqcup \bigcirc$}
               \put(150,10){$)$}
\put(160,10){$= (d^2-1) Y(D), \mbox{ for any graph diagram } D$.}
\end{picture}
\begin{center} {\sc  Figure 2} \end{center}
Throughout this paper,  these  relations will be referred to as the Yamada skein relations.\\

{\bf Theorem 2.1.} {\sl $Y$ is an invariant of ribbon
spatial graphs.}\\
\emph{Proof.} Comparing the recursive formulas defining the Yamada
polynomial $R$ (see \cite{Ya1}, Section 5) to the formulas in
Figure 2, we can easily see that for any spatial graph $\tilde G$
 we have:  $$Y(\tilde
G)(A)=(-d)^{\chi(G)}R(\tilde G)(A^4),$$ where $\chi(G)$ is the
Euler characteristic  of the graph $G$. Thus, $Y$
is an invariant of ribbon spatial graphs. \fin \\
\\
The Kauffman bracket \cite {Ka} polynomial $\prec , \succ$ is an
invariant of regular isotopy of framed links which can be defined
recursively by the following relations:
\begin{center} $\hspace{2cm} \prec \bigcirc \cup
L \succ =d  \prec L \succ, $\\
$\hspace{2cm}  \prec L \succ  =A \prec L_{0} \succ + A^{-1} \prec L_{\infty} \succ ,$\\
\end{center}
where  $d=-A^{2}-A^{-2}$,  $L$, $L_{0}$ and $ L_{\infty}$ are
three links which are identical except in a three-ball where they
are
like in the following picture.\\

\begin{center}
\includegraphics[width=3cm,height=2cm]{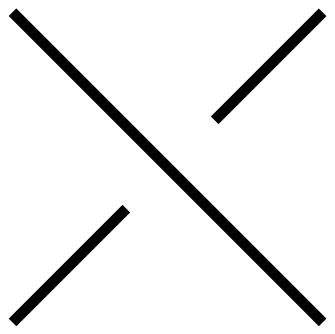} \hspace{1cm}
\includegraphics[width=3cm,height=2cm]{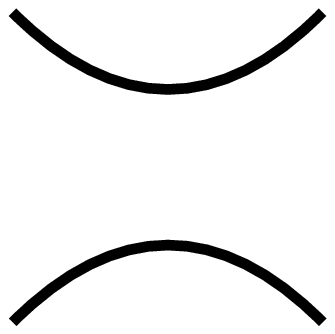} \hspace{1cm}
\includegraphics[width=3cm,height=2cm]{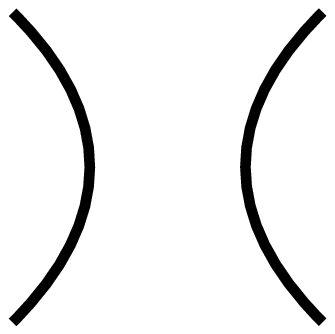}
\end{center}

\begin{center} {\sc Figure 3} \end{center}

It is worth mentioning that the Kauffman bracket is a version of
Jones polynomial for  un-oriented framed links. Murasugi and
Traczyk studied the Jones polynomial of periodic links. They
provided criteria for periodicity of links in terms of the Jones
polynomial.  Here, we write corollaries of these criteria using
the Kauffman bracket,
 see also \cite{Ch1}.\\

{\bf Theorem 2.2 \cite{Mu,Tr}.} {\sl Let $p$ be a prime
and let $L$ be a $p$-periodic link. Then}\\
\begin{tabular}{rl}
(1)& $\prec L \succ \equiv (\prec \overline L \succ)^p$, modulo
$p$,
$d^p - d,$ {\sl where $\overline L$ is the quotient link.} \\
(2)& $\prec L \succ (A) \equiv \prec  L \succ (A^{-1})$, modulo
$p$, $A^{2p} - 1.$
\end{tabular}
\\
\begin{center} {\sc III- Skein modules} \end{center}

Let $M$ be an oriented three-manifold and $\mathcal L$ the set of
all isotopy classes of framed links in $M$. Let $\mathcal
R$=$\Z[A^{\pm 1},d^{-1}]$, where $d=-A^2-A^{-2}$. Let $K(M)$ be
the free $\mathcal R$-module generated by all elements of
$\mathcal L$. We define the Kauffman bracket skein module of $M$,
${\mathcal K}(M)$ to be the quotient of $K(M)$ by the submodule
generated by  all expressions of the form:
\begin{center}
$\hspace{2cm}  \bigcirc \cup
L - d  L $\\
$\hspace{2cm}  L -A L_{0} -A^{-1}  L_{\infty} ,$\\
\end{center}
where  $L$, $L_{0}$ and  $ L_{\infty}$ are three links which are
identical except in a three-ball where they are like in Figure 3.
The Kauffman bracket skein module was computed for several
manifolds. In particular, it was shown that there is an algebra
structure on the skein module of $F \times I$, where $F$ is an
oriented surface. The unit is the empty set and the product is
given by including a copy of $F \times I$ into each of $F\times
[0,1/2]$ and $F \times [1/2,1]$.\\

{\bf Theorem 3.1  \cite{Pr2,Tu}.} {\sl The skein algebra
${\mathcal K}(S^1\times I \times I)$ is isomorphic to the
polynomial algebra ${\mathcal R}[z]$, where $z$ is represented by
a curve in the annulus  as in Figure 4.} \\
\null
\begin{center}
\begin{picture}(0,0)
 \put(0,0){\circle{10}} \put(0,0){\circle{20}} \put(12,0){$z$}
 \put(0,0){\circle{40}}
\end{picture}
\end{center}
\vspace{0.2cm}
\begin{center} {\sc Figure 4} \end{center}
Let $M$ be an oriented three-manifold and $\mathcal G$ the set of
all  isotopy classes of  embeddings of  ribbon graphs in $M$. Let
${\mathcal R}{\mathcal G}$ be the free ${\mathcal R}$-module
generated by ${\mathcal G}$. The Yamada skein module of $M$ which
will be denoted by ${\mathcal Y}(M)$ is defined as the quotient of
the module ${\mathcal R}{\mathcal G}$ by the
submodule generated by  all expressions of the form:\\
\null

               \begin{picture}(0,0)

               \put(135,0){\line(-1,1){13}}
               \put(119,16){\line(-1,1){13}}
               \put(105,0){\line(1,1){30}}
               \put(160,10){$- A^{4}$}
\put(205,15){\oval(28,20)[r]}
               \put(235,15){\oval(28,20)[l]}

               \put(234,10){$ -A^{-4}$}
               \put(300,25){\oval(30,20)[b]}
               \put(300,0){\oval(30,20)[t]}

               \put(324,10){$ +\;d$}
               \put(384,0){\line(-1,1){30}}
               \put(356,0){\line(1,1){30}}

\end{picture}
\\

\begin{picture}(0,0)
                \put(120,10){\line(-1,1){13}}
                \put(110,12){.}
                \put(110,10){.}
                \put(110,8){.}

               \put(120,10){\line(-1,-1){13}}
               \put(120,10){\line(1,0){15}}
               \put(135,10){\line(1,1){13}}
               \put(135,10){\line(1,-1){13}}
               \put(145,12){.}
                \put(145,10){.}
                \put(145,8){.}

               \put(160,10){$-$}
                \put(200,10){\line(-1,1){13}}
               \put(200,10){\line(-1,-1){13}}
               \put(190,12){.}
                \put(190,10){.}
                \put(190,8){.}
               \put(200,10){\line(1,1){13}}
               \put(200,10){\line(1,-1){13}}
               \put(210,12){.}
                \put(210,10){.}
                \put(210,8){.}
                \put(215,10){$\;+\;d^{-1}$}
                \put(275,10){\line(-1,1){13}}
               \put(275,10){\line(-1,-1){13}}
               \put(265,12){.}
                \put(265,10){.}
                \put(265,8){.}
               \put(280,10){\line(1,1){13}}
               \put(280,10){\line(1,-1){13}}
               \put(290,12){.}
                \put(290,10){.}
                \put(290,8){.}

\end{picture}
\\

\begin{picture}(0,0)
\put(120,10){\line(-1,1){13}}
                \put(110,12){.}
                \put(110,10){.}
                \put(110,8){.}
               \put(120,10){\line(-1,-1){13}}
               \put(130,10){\circle{20}}

\put(160,10){$-(d-d^{-1})$}

\put(245,10){\line(-1,1){13}}
                \put(235,12){.}
                \put(235,10){.}
                \put(235,8){.}
               \put(245,10){\line(-1,-1){13}}

\end{picture}
\\

\begin{picture}(0,0)

\put(120,10){$D \bigsqcup \bigcirc$}
 \put(160,10){$- (d^2-1)D , \mbox{ for any ribbon graph  } D.$}
\end{picture}

\begin{center} {\sc  Figure 5} \end{center}

 One can define a graph skein theory for three-manifolds. This theory enjoys the
  same properties  as the skein theory associated to links in three-manifolds,
   see \cite{Pr2}.
 The existence of the Yamada polynomial for spatial graphs  is
equivalent to the fact that ${\mathcal Y}(\S)$ is isomorphic to
$\mathcal R$. As in the case of the Kauffman bracket skein module there is an algebra
structure on the Yamada skein module of $F \times I$, where $F$ is an oriented surface,
 see \cite{Ch2}. \\

{\bf Theorem 3.2.}  {\sl The Yamada skein algebra of the solid
torus ${\mathcal Y}(S^1 \times I \times I)$ is isomorphic to the
polynomial algebra $\cal R$$[z]$, where $z$ is represented by a
curve in the
annulus  as in Figure 4.}\\

{\emph{ Proof.}} In the annulus $S^1 \times I$, we denote by
$b_n$ the bouquet which is made up of  $n$ non-trivial loops, with
the convention that $b_0$ is the empty diagram $\emptyset$. Let
$S_n$ be the graph with 2 vertices and $n+1$ edges as in figure 6.
Finally,  let $\theta_n$ be the $\theta_n$-curve  (contained in some disk).\\
\vspace{1cm}

\begin{picture}(0,0)
               \put(120,0){\oval(80,80)[b]}
               \put(120,0){\oval(80,80)[t]}
                \put(120,0){\circle*{5}}
                \put(120,0){\oval(30,10)[t]}
                \put(120,0){\oval(30,40)[b]}
                \put(120,0){\oval(30,45)[b]}
                \put(120,0){\oval(30,60)[b]}
                \put(120,-25){$.$}
                \put(120,-29){$.$}
                \put(120,-27){$.$}
                \put(120,-55){$S_n$}

                 \put(260,0){\oval(80,80)[b]}
               \put(260,0){\oval(80,80)[t]}
                \put(260,0){\circle*{5}}
                \put(260,-6){\circle{20}}
                \put(260,-8){\circle{22}}
                \put(260,-12){\circle{30}}
                \put(260,-23){$.$}
                \put(260,-27){$.$}
                \put(260,-25){$.$}
                \put(260,-55){$b_n$}

                \put(400,0){\oval(80,80)[b]}
               \put(400,0){\oval(80,80)[t]}
                \put(400,0){\circle*{5}}

                \put(400,0){\oval(30,30)[b]}
                \put(400,0){\oval(30,35)[b]}
                \put(400,0){\oval(30,47)[b]}
                \put(400,-19){$.$}
                \put(400,-21){$.$}
                \put(400,-23){$.$}
                \put(400,-55){$\theta _n$}

\end{picture}

\vspace{2cm}
\begin{center} {\sc Figure 6} \end{center}

 {\bf  Lemma 3.3.} {\sl We have the
  following identities in ${\mathcal Y}(S^1 \times I \times I)$}:
$$\begin{array}{lll}
(i)&\theta_n=&
-d^{-1}\theta_{n-1}+(-\displaystyle\frac{d^2-1}{d})^{n-2}(d^2-1)b_0.\\
(ii)&S_n=&-d^{-1}S_{n-1}+[\displaystyle\frac{-(d^2-1)}{d}]^{n-2}b_1.\\
(iii)& b_n=&S_n+d^{-1}\theta_n.
\end{array}
$$
{ \emph{Proof.}} A routine verification using  the Yamada skein
relations.\fin \\
\\
Let $D$ be a diagram of a spatial graph in the solid torus. It is
easy to see that one can use the first Yamada skein  relation to
transform $D$ onto a linear combination of diagrams   with no
crossings. By using the second Yamada skein relation, it is
possible to transform each of those diagrams onto a linear
combination of diagrams such that each connected component  has
only one vertex. Finally, by applying the third and the fourth
Yamada skein relations, the diagram $D$ may be written as a linear
combination of  bouquets with no contractible cycles. This means
that the set of all bouquets with no contractible cycles,
generates the Yamada skein module of the solid torus.\\
Now, according to (i)  and (ii) of Lemma 3.3, each of $\theta_n$
and $S_n$ may be written as a linear combination  of $b_1$ and
$b_0$. Hence, (iii) implies that $b_n$ can be expressed as a
linear combination of $b_1$ and $b_0$. Thus, we conclude that the
Yamada skein module of the solid torus is generated by all links
in the annulus $S^1\times I$, without trivial components but
including the empty diagram. Finally, the skein algebra of the
solid torus is isomorphic to the polynomial algebra $\cal R$$[z]$.
This ends
the proof of  Theorem 3.2. \fin \\

Throughout the rest of the paper we denote by $\tau _2$ the
Temperly-Lieb algebra with the standard two generators
\begin{picture}(0,0)
\put(10,-10){\line(0,1){25}} \put(17,-10){\line(0,1){25}}
\put(20,0){ and } \put(55,14){\oval(7,17)[b]}
\put(55,-10){\oval(7,17)[t]}
 \put(62,0){.}
\end{picture}
\hspace{23mm} Let
\begin{picture}(0,0)
\put(0,0){$f_{1}=$} \put(30,-10){\line(0,1){25}}
\put(37,-10){\line(0,1){25}} \put(40,0){$-d^{-1}$}
\put(75,14){\oval(7,17)[b]} \put(75,-10){\oval(7,17)[t]}
\put(85,0){ be the Jones-Wenzl projector  in $\tau_2$.}
\end{picture}

\vspace{5mm}
 Let
$G$ be a graph diagram. We define $G'$ to be the linear
combination of link diagrams obtained from $G$ by replacing each
edge of $G$ by two planar strands with a projector $f_1$ in the
cable, and by replacing each vertex of $G$ by a diagram as follows
(the figure illustrates the case
of a four-valent vertex):\\
\begin{picture}(0,0)
\put(184,0){\line(-1,1){30}} \put(156,0){\line(1,1){30}}
\put(220,10){$\longmapsto$} \put(280,20){\framebox(8,4)}
\put(300,20){\framebox(8,4)} \put(280,0){\framebox(8,4)}
\put(300,0){\framebox(8,4)} \put(284,24){\line(-1,1){10}}
\put(303,28){\footnotesize{2}}
 \put(304,24){\line(1,1){10}}
 \put(284,0){\line(-1,-1){10}}
 \put(303,-10){\footnotesize{2}}
 \put(304,0){\line(1,-1){10}}
 \put(295,20){\oval(18,5)[b]}
 \put(282,29){\footnotesize{2}}
 \put(282,-11){\footnotesize{2}}
 \put(295,4){\oval(18,5)[t]}
 \put(283,12){\oval(5,16)[r]}
 \put(307,12){\oval(5,16)[l]}
\end{picture}
\begin{center} {\sc  Figure 7} \end{center}

    Here, writing an integer $n$ beneath an edge $e$ means that this
edge has to be replaced by $n$ parallel ones. \\
 Let $M$ be a three-manifold which is homeomorphic to the product $F \times I$, where $F$ is an
 oriented surface and $I$ is the unit interval.
   Let $\varphi$ be the map from ${\mathcal R}{\mathcal G}$ to $\mathcal
 K$$(M)$ defined on the generators by $\varphi (G)=G'$ and
 extended by linearity to ${\mathcal R}{\mathcal G}$, see also \cite{Ya2} and
 \cite{Yo}. \\
{\bf Lemma 3.4.} {\sl Let $\mathcal Q$ be the sub-module of
 ${\mathcal R}{\mathcal G}$
  generated by the Yamada skein  relations, then $\varphi({\mathcal
 Q})=0$.}\\
 {\emph{Proof.}} For the first skein relation  we refer the reader to
 \cite{KL} page 35. Using the  definition of $f_1$ and the fact that
$f_1^2=f_1$, we should be able to prove the result for the other
skein relations. This is explained by the calculations below:\\
the second skein relation\\
\begin{center}
\includegraphics[width=12cm,height=6cm]{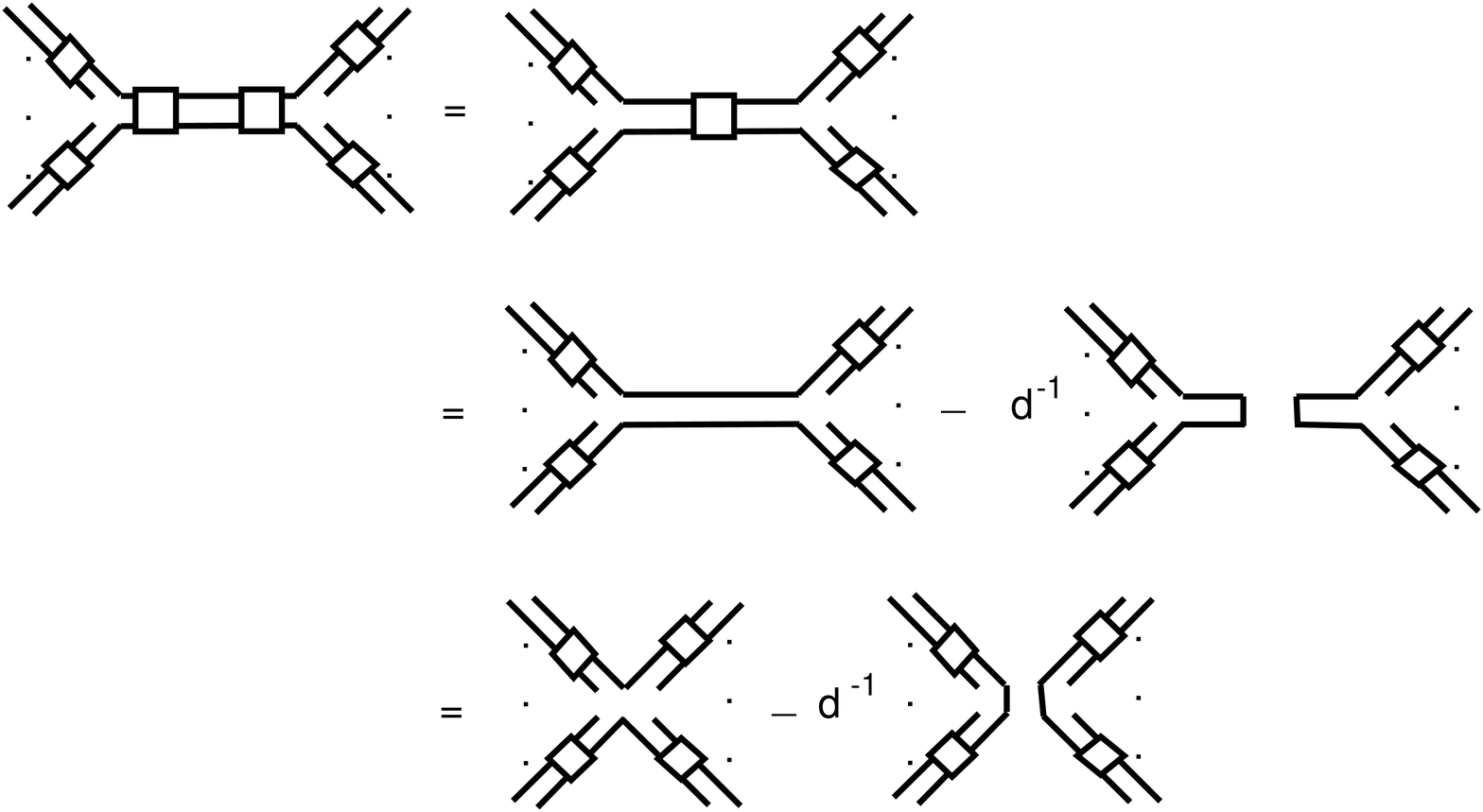}
\end{center}
the third skein relation
\begin{center}
\includegraphics[width=8cm,height=6cm]{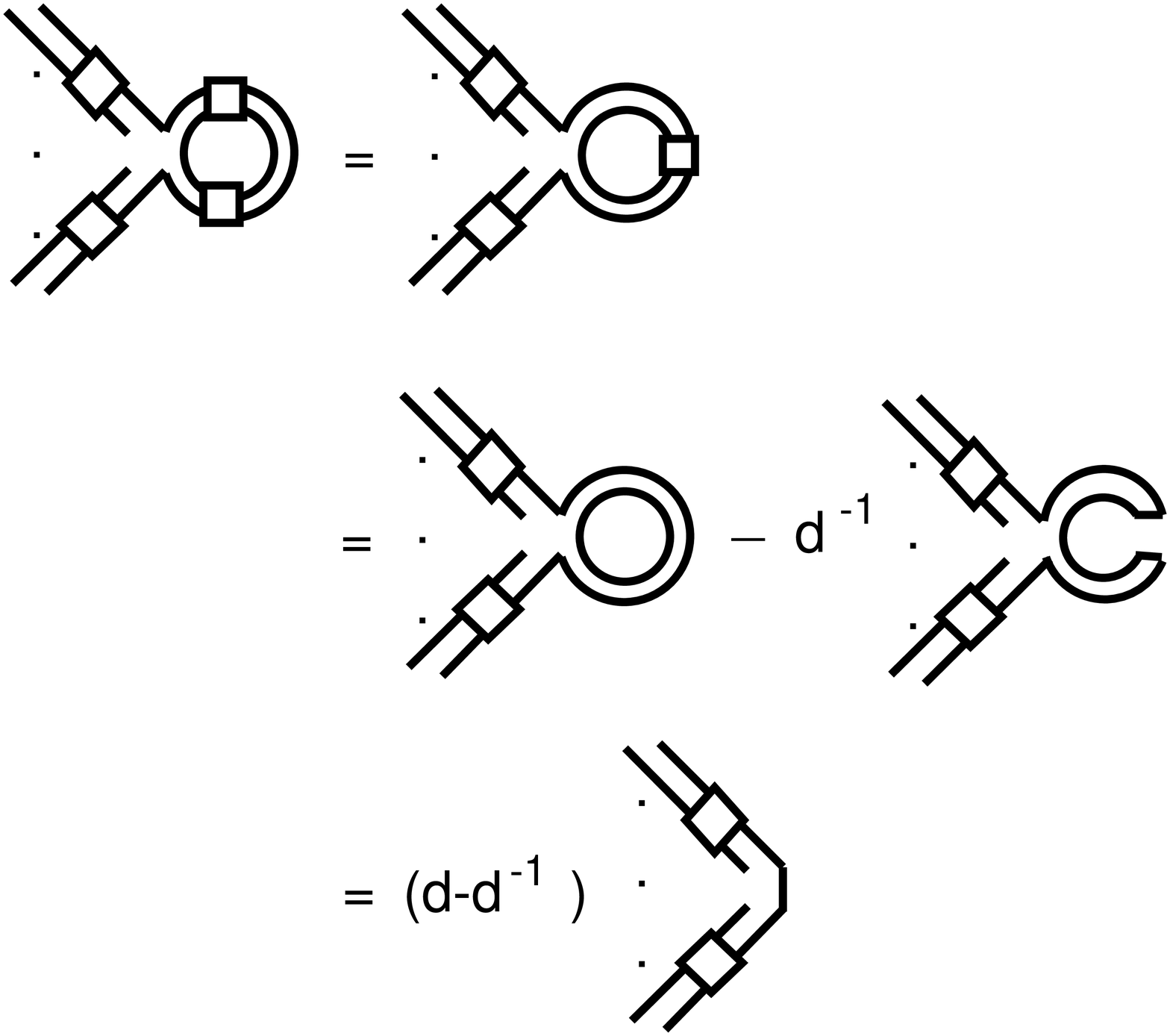}
\end{center}
the fourth skein relation
\begin{center}
\includegraphics[width=6cm,height=3cm]{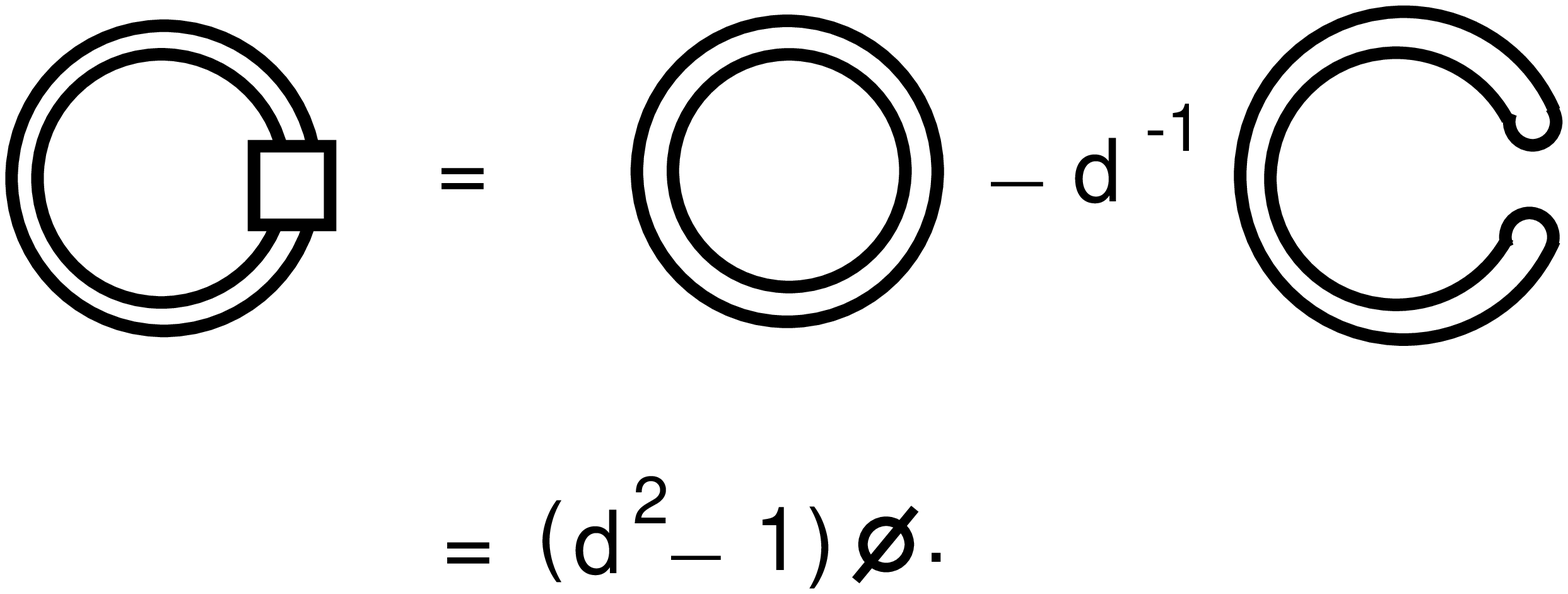}
\end{center}
Since for every ribbon graph $D$ we have   $\varphi(D\cup
\bigcirc)=\varphi(D)\ \varphi(\bigcirc)$, then the identity above
implies that $\varphi(D\cup \bigcirc - (d^2-1)D)=0$. This
ends the proof of Lemma 3.4. \fin\\

By Lemma 3.4, $\varphi$ defines a map from the Yamada
 skein module  $\mathcal
 Y$$(S^1 \times I \times I)$ to the Kauffman bracket skein module $\mathcal
 K$$(S^1 \times I \times I)$, we will denote this map  by $\Phi$.
 It is easy to see that $\Phi$ is a homomorphism of algebras. Moreover, we have\\

{\bf Lemma 3.5.} {\sl The algebra homomorphism  $\Phi$ is injective.}\\
{\emph{ Proof.}} Direct computation shows that  $\Phi (z)=z^2-1.$
Hence, $\Phi$ is injective.\fin \\

\begin{center} {\sc IV- Proof of Theorem 1.1} \end{center}
Let $\tilde G$ be a $p$-periodic spatial graph. Let $L_{\tilde G}$
be the set of all links that appear when we expand $\tilde G$ as a
linear combination of links using the map $\Phi$. Namely:
$$\Phi (\tilde G)= \displaystyle\sum_{D \in L_{\tilde G}}{}
d^{\alpha_D}D,$$ where $\alpha_{D}$ is an integer.\\
As there is an action of the finite cyclic group of order $p$ on
the set of vertices and on the set of crossings  of $\tilde G$, we
can see easily that this group acts on $L_{\tilde G}$. Moreover,
there are two
types of orbits:\\
Orbits made up of only one link which is  $p$-periodic. Let
$L_{p,\tilde G}$ be the subset of $L_{\tilde G}$ made
 up of $p$-periodic
links.\\
Orbits made up of $p$ links which are cyclically permuted by the
rotation.\\
In the second case the $p$ links are the same. Moreover, if $D$
and $D'$ belong to the same orbit then $\alpha_D = \alpha_{D'}$.
Hence  the contribution of the links from the same orbit adds to
zero modulo $p$. Consequently, we have the following congruence
modulo $p$.
$$\Phi (\tilde G)\equiv \displaystyle\sum_{D \in L_{p,\tilde G}}{}
d^{\alpha_D}D.$$

If $D$ is $p$-periodic then by Theorem 2.2, we have the following
identity in the Kauffman bracket skein module
$$D \equiv \overline D ^p, \mbox{ modulo } p, d^p -d.$$
Thus we have the following congruence modulo $p$ and $d^p -d$:
$$\Phi (\tilde G)\equiv \displaystyle\sum_{D \in L_{p,\tilde G}}{}
d^{\alpha_D} \overline D ^p.$$

It is easy to see that if ${D \in L_{p,\tilde G}}$ then $\alpha_D$
can be written $p \alpha'_D$ for some integer  $\alpha'_D$. If we
expand the graph diagram $\underline { \tilde G}$ into a linear
combination of links using the map $\Phi$, we can see that the
links that appear in the sum are exactly the quotients of the
elements of $L_{p,\tilde G}$. Moreover, if $D$ appear in $\Phi
(\tilde G)$ with coefficient $\alpha_D$, then $\overline D$
appears in $\Phi (\underline {\tilde G})$, with coefficient
$\alpha'_D$. This allows us to conclude that:
$$\Phi (\tilde G) \equiv (\Phi (\underline {\tilde G}))^p,
 \mbox{ modulo } p, d^p -d.  $$
Using the fact that $\Phi$ is an injective  homomorphism of
algebras we conclude that in the Yamada skein module of the solid
torus we
have:\\
$$\tilde G\equiv (\underline {\tilde G})^p \mbox{ modulo } p, d^p -d.  $$
This ends the proof of part (1) of Theorem 1.1\\
Similar arguments are used to prove the second part of Theorem 1.1.\\

\end{document}